\documentclass[10pt,journal]{IEEEtran}
\usepackage{graphicx}          
\usepackage{amsmath}           
\usepackage{amssymb}
\usepackage{mathrsfs}
\usepackage{algorithm}
\usepackage{algpseudocode}
\usepackage{hyperref}
\usepackage[latin1]{inputenc}  

\usepackage{amsthm}

\theoremstyle{definition}

\theoremstyle{remark}

\theoremstyle{plain}

\DeclareMathOperator{\E}{E}
\DeclareMathOperator*{\argmin}{arg\,min}

\newcommand{\mbf}[1]{\mathbf{#1}}
\newcommand{\mbs}[1]{\boldsymbol{#1}}
\newcommand{\what}[1]{\widehat{#1}}

\begin{document}

\title{Line Spectrum Estimation with Probabilistic Priors}
\author{Dave Zachariah, Petter Wirfält, Magnus Jansson and Saikat Chatterjee\thanks{The authors are with the ACCESS Linnaeus Centre, KTH Royal Institute of Technology, Stockholm. E-mail: $\{$dave.zachariah, wirfalt, magnus.jansson$\}$@ee.kth.se and saikatchatt@gmail.com. This work was partially supported by the Swedish Research Council under contract 621-2011-5847. The research leading to these results has received funding from the European Research Council under the European Community's Seventh Framework Programme (FP7/2007-2013) / ERC grant agreement n$^\circ$ 228044. }}

\maketitle

\begin{abstract}
For line spectrum estimation, we derive the maximum a posteriori probability estimator where prior knowledge of frequencies is modeled probabilistically. Since the spectrum is periodic, an appropriate distribution is the circular von Mises distribution that can parameterize the entire range of prior certainty of the frequencies. An efficient alternating projections method is used to solve the resulting optimization problem. The estimator is evaluated numerically and compared with other estimators and the Cramér-Rao bound.
\end{abstract}


\section{Introduction}

Line spectrum estimation is a classical problem in signal processing
with many applications, including communications, radar, sonar and
seismology \cite{Stoica&Moses2005}. In such applications, the observed
signal contains frequencies that may be known to varying degrees of
certainty. Examples include diagnosis applications where the power line frequency may appear \cite{WirfaltEtAl2011}; communications systems with known carrier frequencies; characterization of circuits, e.g., analogue to digital converters, power amplifiers, etc., using sinusoidal test signals which further give rise to known harmonics. Assuming a subspace approach, the problem is similar to
direction of arrival estimation with uniform linear arrays, for which
methods that incorporate prior knowledge have been developed \cite{LinebargerEtAl1995,BouleuxEtAl2009,SteinwandtEtAl2011}.
In \cite{WirfaltEtAl2011}, the approach of \cite{BouleuxEtAl2009} was developed for frequency estimation. This method, however, assumes perfect, deterministic knowledge of a subset of frequencies while assuming no prior knowledge about the remaining ones. Alternative methods for incorporating prior knowledge in line spectrum estimation have been developed by imposing sparsifying penalties on spectral amplitudes over a grid of frequencies, cf. \cite{BourguignonEtAl2007}.

In this paper, we approach the problem in a probabilistic manner, in
which the prior certainty of each frequency can vary. For discrete-time signals, any inferred value of the dimensionless frequency $\omega$ is equivalent to that of $\omega' = \omega + 2 \pi k$, where $k$ is an integer. Thus for consistent probabilistic inference, the prior and posterior probability density functions (pdf) need to be periodic or circular \cite{Lee2010}. Probabilistic treatment of the frequencies is uncommon in the literature. A rare example is \cite{Bretthorst1988} but it assumes  noninformative priors; in that case, for well-separated frequencies, the resulting maximum a posteriori probability (MAP) estimator is the periodogram. See also \cite{Dou&Hodgson1995}.

We derive the MAP line spectrum estimator that exploits prior knowledge of the frequencies. This information may be given from past experience or in the process of detecting the number of cisoids. The information is then particularly useful when few samples are available and/or when the signal-to-noise ratio is low, i.e., in conditions where standard line spectrum estimators may fail. A computationally efficient alternating projections method is used to solve the resulting optimization problem. The performance of the MAP estimator is evaluated numerically and compared with two other estimators, the Cramér-Rao bound (CRB) and the hybrid CRB.

\emph{Notation:} $\mbf{A}^*$ and $\mbf{A}^\dagger$ denote the Hermitian transpose and Moore-Penrose pseudo-inverse of the matrix $\mbf{A}$, respectively. $\mbs{\Pi}_{\mbf{A}}$ and
$\mbs{\Pi}^\perp_{\mbf{A}}$ denote the orthogonal projection matrices
onto the range space of $\mbf{A}$ and its complement,
respectively. $\text{tr}\{ \cdot \}$ denotes the trace operator. $\mbf{e}_i$ is the $i$th standard basis vector in $\mathbb{R}^m$.

\section{Problem formulation}

A set of $m$ samples of a sum of $d$ cisoids is observed
\begin{equation}
y(t) = \sum^d_{i=1} s_i e^{j \omega_i t} + n(t) \in \mathbb{C}, \quad
t = 0, \dots, m-1,
\label{eq:scalarmodel}
\end{equation}
where $s_i \in \mathbb{C}$ parameterizes the amplitude and phase of
the $i$th cisoid, and $\omega_i$ is its frequency. The zero-mean noise $n(t)$ is assumed to be independent and identically distributed (i.i.d.) complex Gaussian with variance $\sigma^2$. For identifiability $m > d$. In vector form, \eqref{eq:scalarmodel} can
be written as
\begin{equation}
\mbf{y} = \mbf{A}(\mbs{\omega}) \mbf{s} + \mbf{n} \in \mathbb{C}^m
\end{equation}
where $\mbf{n} = [n(0) \: \cdots \: n(m-1)]^\top \in \mathbb{C}^m$ and
$\mbf{s} = [s_1 \: \cdots \: s_d]^\top \in \mathbb{C}^d$.
The Vandermonde matrix
\begin{equation*}
\mbf{A}(\mbs{\omega}) = \begin{bmatrix}
1 & \cdots & 1 \\
e^{j \omega_1} & \dots & e^{j \omega_d} \\
\vdots & \ddots & \vdots \\
e^{j (m-1) \omega_1} & \cdots & e^{j (m-1)  \omega_d}
\end{bmatrix}
\in \mathbb{C}^{m \times d}
\end{equation*}
is parameterized by $\mbs{\omega} = [ \omega_1 \cdots
\omega_d]^\top$. The goal is to estimate $\mbs{\omega}$,~$\mbf{s}$ and $\sigma^2$ from $\mbf{y}$.

No prior knowledge of $\mbf{s}$ or $\sigma^2$ is assumed. We model this using noninformative Jeffreys priors $p(\mbf{s}) \propto 1$ and $p(\sigma^2) \propto 1/\sigma^2$ \cite{Tiao&Zellner1964}, which enables a consistent Bayesian treatment of the estimation problem.

\begin{figure}
  \begin{center}
    \includegraphics[width=1.0\columnwidth]{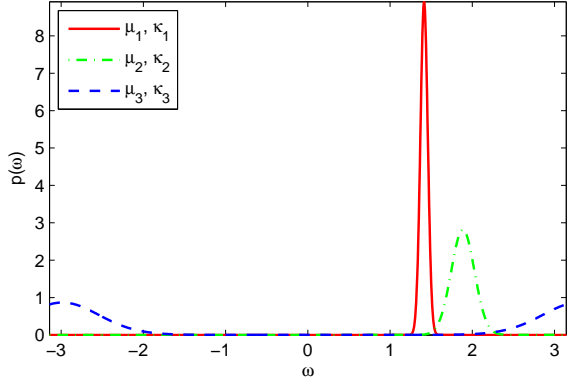}
  \end{center}
  \caption{Illustration of three different prior pdfs over frequencies
    $\omega \in [-\pi, \pi)$. The circular means are $\mu_1 = 0.45\pi$, $\mu_2 =  0.60 \pi$ and $\mu_3 = -0.95\pi$. The dispersion is parameterized by $\kappa_i$ and here set to $\kappa_1 = 500, \kappa_2 = 50$ and $\kappa_3 = 5$. For further illustrations of the von Mises pdf, see \cite{AbdiEtAl2002}.}
  \label{fig:vonMises}
\end{figure}

The frequencies $\{ \omega_i \}$ are modeled as independent random variables,
with circular pdfs, such that $p(  \omega_i ) = p( \omega_i + 2\pi k
)$ for any integer $k$. A tractable prior distribution with this
property is the von Mises distribution, $\omega_{i} \sim
\mathcal{M}(\mu_{i}, \kappa_{i})$, which can be thought of as a
circular analogue of the Gaussian distribution on the line
\cite{Lee2010}. Its pdf is
\begin{equation}
p(\omega_i; \mu_i, \kappa_i) = \frac{1}{2\pi I_0(\kappa_i)} e^{\kappa_i \cos(\omega_i - \mu_i)},
\end{equation}
with $I_\ell(\kappa_i)$ being the modified Bessel function of order $\ell$.
The circular mean and circular variance of $\omega_i$ are $\E[\omega_i]=\mu_i$
and $\text{Var}[\omega_i] = 1 - I_1(\kappa_i) / I_0(\kappa_i)$, respectively.\footnote{For circular distributions the $n$th trigonometric moment is defined by $\varsigma_n = \E_{\theta}[e^{jn \theta}]$. In polar coordinates, $\varsigma_1 = \rho e^{j \mu}$, where $\mu$ and $1-\rho$ define the circular mean and variance of $\theta$, respectively \cite{Lee2010}.}
As the concentration parameter is varied to its extremes, $\kappa_i \rightarrow 0$ and $\kappa_i \rightarrow \infty$, the pdf of $\omega_i \in [-\pi, \pi)$ approaches a uniform and a Gaussian pdf with variance $1/\kappa_i$, respectively \cite{EvansEtAl2000}. Thus, $\kappa_i$ enables parametrization of the prior certainty of frequency $\omega_i$, from complete ignorance to virtual certainty. An illustration of the von Mises pdf is given in Fig.~\ref{fig:vonMises}. For relatively small variances, the approximate properties of the pdf provides a practical way to select $\kappa_i$ using the confidence level of a Gaussian with variance $1/\kappa_i$.

\section{MAP estimator}
The joint maximum a posteriori probability estimator of $\mbs{\omega}$,~$\mbf{s}$ and $\sigma^2$ is given by maximization of $p(\mbs{\omega}, \mbf{s}, \sigma^2| \mbf{y})$, or equivalently the cost function $J( \mbs{\omega}, \mbf{s}, \sigma^2) = J_1( \mbs{\omega}, \mbf{s}, \sigma^2) + J_2( \mbs{\omega})$, where
\begin{equation}
\begin{split}
J_1(\mbs{\omega}, \mbf{s},\sigma^2) = \ln p(  \mbf{y} | \mbs{\omega}, \mbf{s}, \sigma^2 )  + \ln p(\mbf{s}) + \ln p(\sigma^2)
\end{split}
\end{equation}
and $J_2(\mbs{\omega}) = \ln p(\mbs{\omega})$. Using the noninformative priors of $\mbf{s}$
and $\sigma^2$, one obtains
\begin{equation}
J_1(\mbs{\omega}, \mbf{s},\sigma^2) = -m \ln \sigma^2 - \frac{1}{\sigma^2} \| \mbf{y} - \mbf{A}(\mbs{\omega})\mbf{s} \|^2_2 - \ln \sigma^2 + K_1,
\end{equation}
where $K_1$ is a constant and the maximizers are given by $\hat{\mbf{s}} = \mbf{A}^\dagger (\mbs{\omega}) \mbf{y}$ and $\hat{\sigma}^2 = \mbf{y}^* \mbs{\Pi}^\perp_{\mbf{A}}(\mbs{\omega}) \mbf{y}/(m+1)$. Then
\begin{equation}
J_1(\mbs{\omega}, \hat{\mbf{s}},\hat{\sigma}^2) = -(m+1) \ln \left(  \mbf{y}^* \mbs{\Pi}^\perp_{\mbf{A}}(\mbs{\omega}) \mbf{y} \right) + K'_1,
\label{eq:J1}
\end{equation}
where $K'_1 = K_1 - (m+1)$. Further, as the frequencies are independently distributed, $J_2(\mbs{\omega})$ can be
written compactly in terms of $\mbf{A}(\mbs{\omega})$,
\begin{equation}
\begin{split}
J_2(\mbs{\omega}) &= \sum^{d}_{i=1} \ln p( \omega_i ) \\
&= \sum_i \kappa_i \cos(\omega_i - \mu_i)  + K_2 \\
&= \sum_i \frac{\kappa_i}{2} e^{j(\omega_i - \mu_i)} + \frac{\kappa_i}{2} e^{-j(\omega_i - \mu_i)} + K_2 \\
&= \text{Re} \{ \mbf{e}^*_2 \mbf{A}(\mbs{\omega}) \mbs{\pi} \} + K_2,
\end{split}
\label{eq:J2}
\end{equation}
where $\mbs{\pi} = [ \pi_1 \cdots \pi_d]^\top$ with $\pi_i = \kappa_i e^{-j \mu_i}$ parameterizing the prior knowledge, and $K_2$ is a constant.

\subsection{Concentrated cost function}
Combining \eqref{eq:J1} and \eqref{eq:J2}, the MAP estimator is given by
\begin{equation}
\what{\mbs{\omega}}_{\text{map}} = \argmin_{\mbs{\omega} \in \Omega}
V_{\text{map}}(\mbs{\omega}),
\label{eq:J_map}
\end{equation}
where
\begin{equation*}
\begin{split}
V_{\text{map}}(\mbs{\omega}) &\triangleq  \left(  \mbf{y}^*
  \mbs{\Pi}^\perp_{\mbf{A}} (\mbs{\omega}) \mbf{y} \right) e^{\phi(\mbs{\omega})},
\end{split}
\end{equation*}
$\phi(\mbs{\omega}) = \text{Re}\{ \mbf{e}^*_2 \mbf{A}(\mbs{\omega}) \mbs{\beta} \} $
and $\mbs{\beta} \triangleq -\frac{1}{m+1}[ \kappa_1 e^{-j \mu_1}
\cdots \kappa_d e^{-j \mu_d}]^\top$. The cost function is highly nonlinear and multimodal, but given a good initial guess the optimization problem can be solved by a grid or Newton-based search method \cite{Stoica&Moses2005}. The computational complexity of such a $d$-dimensional optimization problem may, however, be prohibitive. For this reason we formulate an alternating projection method based on \cite{Ziskind&Wax1988}, which reduces the problem to a series of 1-dimensional optimization problems and results in a computationally tractable estimator.

\subsection{Alternating projection solution}

Let the $i$th column of $\mbf{A}$ be denoted as $\mbf{a}_i \in \mathbb{C}^{m \times 1}$, corresponding to $\omega_i$, and the remaining columns $\mbf{A}_i$, corresponding the remaining frequencies denoted $\check{\mbs{\omega}}_i$. Then the projection operator can be decomposed as
$\mbs{\Pi}_{\mbf{A}_i : \mbf{a}_i} = \mbs{\Pi}_{\mbf{A}_i} +
\mbs{\Pi}_{\tilde{\mbf{a}}_i}$, where $\tilde{\mbf{a}}_i =
\mbs{\Pi}^\perp_{\mbf{A}_i} \mbf{a}_i$, so that
$\mbs{\Pi}_{\tilde{\mbf{a}}_i} = \tilde{\mbf{a}}_i \tilde{\mbf{a}}^*_i
/ \| \tilde{\mbf{a}}_i \|^2$.

The cost function $V_{\text{map}}(\mbs{\omega})$ is minimized for each frequency $\omega_i$, holding $\check{\mbs{\omega}}_i$ constant. Hence the $d$-dimensional optimization problem \eqref{eq:J_map} is relaxed into an iteration of 1-dimensional grid searches:
\begin{equation}
\widehat{\omega}_i = \argmin_{\omega_i \in \Omega_i } V(\omega_i; \check{\mbs{\omega}}_i ),
\label{eq:search}
\end{equation}
where
\begin{equation}
\begin{split}
V(\omega_i; \check{\mbs{\omega}}_i ) &\triangleq \left( \mbf{y}^* \mbs{\Pi}^\perp_{\mbf{A}_i}\mbf{y} - \frac{| \mbf{y}^* \mbs{\Pi}^\perp_{\mbf{A}_i} \mbf{a}(\omega_i) |^2}{ \| \mbs{\Pi}^\perp_{\mbf{A}_i} \mbf{a}(\omega_i)\|^2 } \right) e^{\phi_i(\omega_i)},
\end{split}
\label{eq:J_ap}
\end{equation}
$\mbf{a}(\omega) = \begin{bmatrix} 1 & e^{j\omega} & \cdots & e^{j(m-1)\omega} \end{bmatrix}^\top$ and $\phi_i(\omega)= \text{Re}\{ \beta_i e^{j \omega} \}$. This follows from the decomposition of $\mbs{\Pi}_{\mbf{A}}(\mbs{\omega})$ and $\text{Re}\{ \mbf{e}^*_2 \mbf{A}(\mbs{\omega}) \mbs{\beta}\} = \text{Re}\{\mbf{e}^*_2 \beta_i \mbf{a}_i(\omega_i)\} + \text{Re}\{\mbf{e}^*_2 \sum_{\ell \neq i} \beta_\ell \mbf{a}_\ell(\omega_\ell) \}$ in \eqref{eq:J_map}. Note that, as $\check{\mbs{\omega}}_i$ is held constant, the latter term is removed from the cost function.

The search \eqref{eq:search} is performed sequentially for all $i = 1, \dots, d$ over a grid of $g$ points, denoted $\Omega_i$. The grid searches are repeated until the difference between iterates, $|\Delta \widehat{\omega}_i|$, is less than some $\varepsilon$. A key element in the algorithm is the initialization, and we follow the procedure of \cite{Ziskind&Wax1988}. To reduce the initial error in the search incurred when holding $\check{\mbs{\omega}}_i$ constant, the algorithm is initialized by setting $\what{\mbs{\omega}} = \varnothing$ and with frequencies $i = 1, \dots, d$ sorted in descending order with respect to their prior certainty, as quantified by the magnitude of $\beta_i$. Then the estimates are initialized sequentially $i = 1,\dots, d$,
\begin{equation*}
\what{\omega}_i = \argmin_{\omega_i \in \Omega_i } \; V(\omega_i; \what{\mbs{\omega}} ) \; \text{followed by} \; \what{\mbs{\omega}} := \what{\mbs{\omega}} \cup \what{\omega}_i.
\end{equation*}
The grids $\Omega_i$ are initially in the interval $[-\pi, \pi)$ and can subsequently be refined by narrowing the intervals centered around the previous estimates, $\widehat{\omega}_i$. The refinement is repeated $L$ times. The alternating projections-based MAP estimator is summarized in Algorithm~\ref{alg:MAP}.

\begin{algorithm}
\caption{MAP line spectrum estimator} \label{alg:MAP}
\begin{algorithmic}[1]
\State Input: $\mbf{y}, \{ \mu_i, \kappa_i \}^d_{i=1}$ and $L$
\State Initialize $\what{\mbs{\omega}} = \varnothing$ and form $\mbs{\beta}$ and $\Omega^1_i$
\For{$\ell = 1,\dots, L$}
    \Repeat
        \State For $i = 1, \dots, d$
        \State Form $\mbs{\Pi}^\perp_{\mbf{A}_i}(\check{\mbs{\omega}}_i)$
        \State $\widehat{\omega}_i = \argmin_{\omega_i \in
          \Omega^\ell_i } V(\omega_i ;
        \check{\mbs{\omega}}_i )$ using \eqref{eq:J_ap}
    \Until{ $|\Delta \widehat{\omega}_i| < \varepsilon_\ell$ }
    \State Refine $\Omega^{\ell+1}_i, \forall i$
\EndFor
\State $\hat{\mbf{s}} = \mbf{A}^\dagger(\widehat{\mbs{\omega}}) \mbf{y}$
\State $\hat{\sigma}^2 = \mbf{y}^* \mbs{\Pi}^\perp_{\mbf{A}}(\widehat{\mbs{\omega}}) \mbf{y}/(m+1)$
\State Output: $\widehat{\mbs{\omega}}$, $\hat{\mbf{s}}$ and $\hat{\sigma}^2$
\end{algorithmic}
\end{algorithm}

\section{Experimental results}

The estimator is evaluated by means of simulation with respect to the root mean square error, $\text{RMSE}(\what{\omega}_i) \triangleq \sqrt{ \E[ \tilde{\omega}^2_i ] }$, where $\tilde{\omega}_i$ is the estimation error (modulo-$2 \pi$) for a given realization of $\mbf{n}$, $\mbf{s}$ and $\mbs{\omega}$. The RMSE is estimated using $10^4$ Monte Carlo runs.

The Cramér-Rao bound (CRB) for deterministic $\mbs{\omega}$ and $\mbf{s}$ is given by
$\mbf{C}(\mbs{\omega}, \mbf{s}) =  \left( \frac{2}{\sigma^2}
  \text{Re}\left \{ \mbf{S}^* \mbf{D}^* \mbf{\Pi}^\perp_{\mbf{A}}
    \mbf{D}\mbf{S}  \right \} \right)^{-1},$ where $\mbf{S} = \text{diag}(\mbf{s})$ and the $i$th column of $\mbf{D}$ is $\mbf{d}_i
= \frac{d\mbf{a}(\omega_i)}{d\omega_i}$
\cite{Stoica&Nehorai1989}. In the simulations, $\mbf{C}(\mbs{\omega}, \mbf{s})$ is
averaged over all realizations of $\mbs{\omega}$ and $\mbf{s}$. For conditionally
unbiased estimators, the diagonal elements set the limit
$\text{RMSE}(\what{\omega}_i) \geq \sqrt{c_{ii}}$. The posterior Cramér-Rao bound does not exist for stochastic frequencies since regularity conditions do not hold for the von Mises pdf \cite{VanTrees2001, Routtenberg&Tabrikian2012}. When $\kappa_i$ are large, however, the variances of the frequencies are small and further $p(\omega_i)$ can be approximated locally by a Gaussian with variance $1/\kappa_i$. Then, following \cite{WahlbergEtAl1991}, we can formulate an approximate hybrid Cramér-Rao bound (ACRB) \cite{Rockah&Schultheiss1987,Reuven&Messer1997}, $\mbf{C}(\mbf{s}) \simeq \left( \frac{2}{\sigma^2} \text{Re}\left \{ \mbf{S}^* \bar{\mbf{D}}^* \bar{\mbf{\Pi}}^\perp_{\mbf{A}} \bar{\mbf{D}}\mbf{S}  \right \} + \mbs{\Lambda}_\omega \right)^{-1}$, where $\bar{\mbf{D}}$ and $\bar{\mbf{\Pi}}^\perp_{\mbf{A}}$ are evaluated at the mean frequencies and $\mbs{\Lambda}_\omega = \text{diag}\{ \lambda_{\omega,1}, \dots , \lambda_{\omega,d} \}$ in which $\lambda_{\omega,i}$ equals $\kappa_i$ or 0 depending on whether the frequency $\omega_i$ is treated stochastically or deterministically, respectively.

For further comparison we also consider the `Estimation of Signal Parameters via Rotational Invariance Techniques' (\textsc{Esprit}) estimator, using the forward-backward covariance estimate \cite{RoyEtAl1986,Stoica&Moses2005}, which does not take prior knowledge into account, and the Markov-based \textsc{Pledge} estimator, which is state of the art for deterministic prior knowledge \cite{ErikssonEtAl1994, WirfaltEtAl2011}.

\subsection{Setup}

We consider $d=3$ cisoids with varying prior knowledge of the frequencies. The prior certainties of $\omega_1$, $\omega_2$ and $\omega_3$ are parameterized by concentration parameters $\kappa_1 = 2\cdot10^3$, $\kappa_2 = 2\cdot10^2$ and $\kappa_3 = 0$; corresponding to  standard deviations of approximately $7\cdot 10^{-3}\pi$ and $2\cdot 10^{-2} \pi$ radians, and complete ignorance, respectively. Frequencies $\omega_1$ and $\omega_2$ are randomly generated with circular means $\mu_1 = 0.45\pi$ and $\mu_2 = 0.60 \pi$ \cite{Berens2009}, while $\omega_3$ is set deterministically to $0.75\pi$ (as the estimator is ignorant, $\kappa_3=0$, it can set $\mu_3$ arbitrarily, e.g., $\mu_3 = 0$.). This choice prevents realizations of randomly generated frequency separations well below the resolution limit of the periodogram resulting in near-degeneracy of the estimation problem.

The cisoid amplitude and phase were set as $s_i = \alpha_i e^{j \varphi_i}$, where $\alpha_i \equiv 1$ and $\varphi_i$ is drawn uniformly over $[0, 2\pi)$ for each realization. Two signal parameters are varied: (i) the signal-to-noise ratio, $\text{SNR}_i \triangleq \E[|s_i|^2] / \E[|n(t)|^2] = \sigma^{-2}$ and (ii) the number of samples, $m$.

For the MAP estimator, a grid of $g=500$ points was used. The algorithm was set to terminate after $L = 10$ refinement levels. The refinement level is performed by reducing the search segment by half, resulting in a resolution limit of about $\pi/(2^{L-1}g) \approx 4 \cdot 10^{-6} \pi$ radians. In our experience, a convergence tolerance $\varepsilon$ of 2 grid points prevents occasional cycling of the minimum point of \eqref{eq:J_ap}. For \textsc{Pledge} we set $\mu_1$ as the prior of $\omega_1$. For \textsc{Esprit} and \textsc{Pledge}, the window length was fixed at $m/2$.

\subsection{Results}

\begin{figure}
  \begin{center}
    \includegraphics[width=1.00\columnwidth]{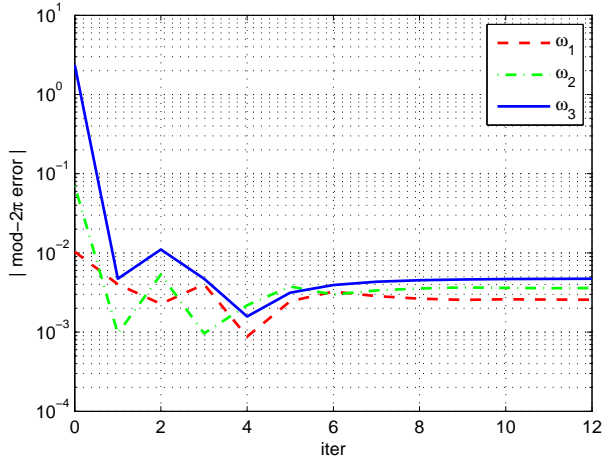}
  \end{center}
  \vspace{-0.5cm}
  \caption{The convergence of the MAP estimates. Absolute error $|\tilde{\omega}_i |$ versus iteration for a typical realization with $m=32$ and SNR = 10~dB. The algorithm terminated at the 12th iteration.}
  \label{fig:convergence}
\end{figure}

An illustration of the convergence of the MAP estimator is given in Fig.~\ref{fig:convergence}. Under the same setting, $m=32$ samples are processed in approximately 4, 200 and 740 milliseconds using the current implementations of \textsc{Esprit}, \textsc{Pledge} and MAP, respectively.

\begin{figure}
  \begin{center}
    \includegraphics[width=1.0\columnwidth]{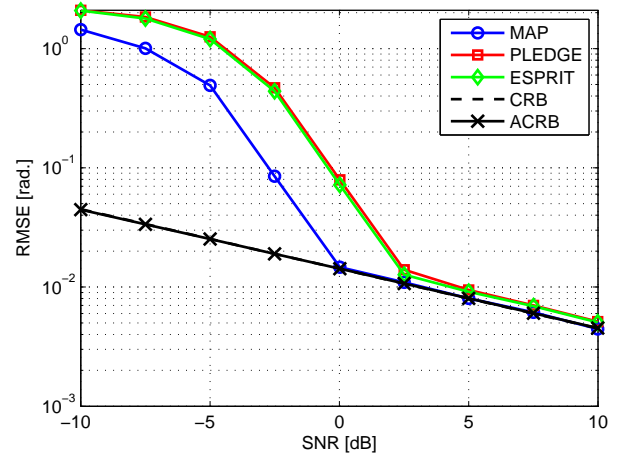}
  \end{center}
  \vspace{-0.5cm}
  \caption{RMSE$(\what{\omega}_3)$ vs SNR for $m=32$.}
  \label{fig:SNR_03_m32}
\end{figure}

The results from the Monte Carlo runs are given below. First, Figs.~\ref{fig:SNR_01_m32}, \ref{fig:SNR_02_m32}, and \ref{fig:SNR_03_m32} show the RMSE of the frequency for each cisoid when fixing $m=32$ and varying SNR. Recall that the prior certainty of the frequencies is in decreasing order. The MAP estimator incorporates the information optimally and is therefore capable of producing estimates of $\omega_1$ that converge to the CRB from below as SNR increases in Fig.~\ref{fig:SNR_01_m32}. The intermediate case, $\omega_2$, is illustrated in Fig.~\ref{fig:SNR_02_m32}. For the frequency with minimum certainty, $\omega_3$, MAP closes the gap to the bound faster than \textsc{Esprit}.

In this scenario the random frequency separation is on average wide so that the average performance improvements of \textsc{Pledge} over \textsc{Esprit} are marginal at low SNR, cf. deterministic scenario in \cite{WirfaltEtAl2011}. For the first cisoid, \textsc{Pledge} cannot improve on the prior of $\omega_1$, which it takes to be perfect deterministic knowledge. Hence the actual deviations of $\omega_1$ from $\mu_1$ impedes the estimates of $\omega_2$ and $\omega_3$ at high SNR.

\begin{figure}
  \begin{center}
    \includegraphics[width=1.0\columnwidth]{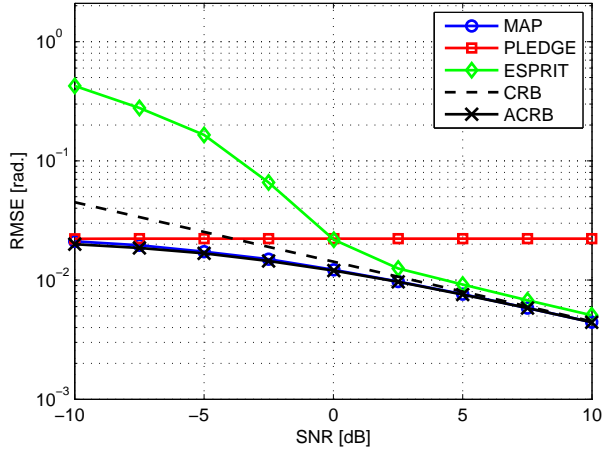}
  \end{center}
  \vspace{-0.5cm}
  \caption{RMSE$(\what{\omega}_1)$ vs SNR for $m=32$.}
  \label{fig:SNR_01_m32}
\end{figure}

\begin{figure}
  \begin{center}
    \includegraphics[width=1.0\columnwidth]{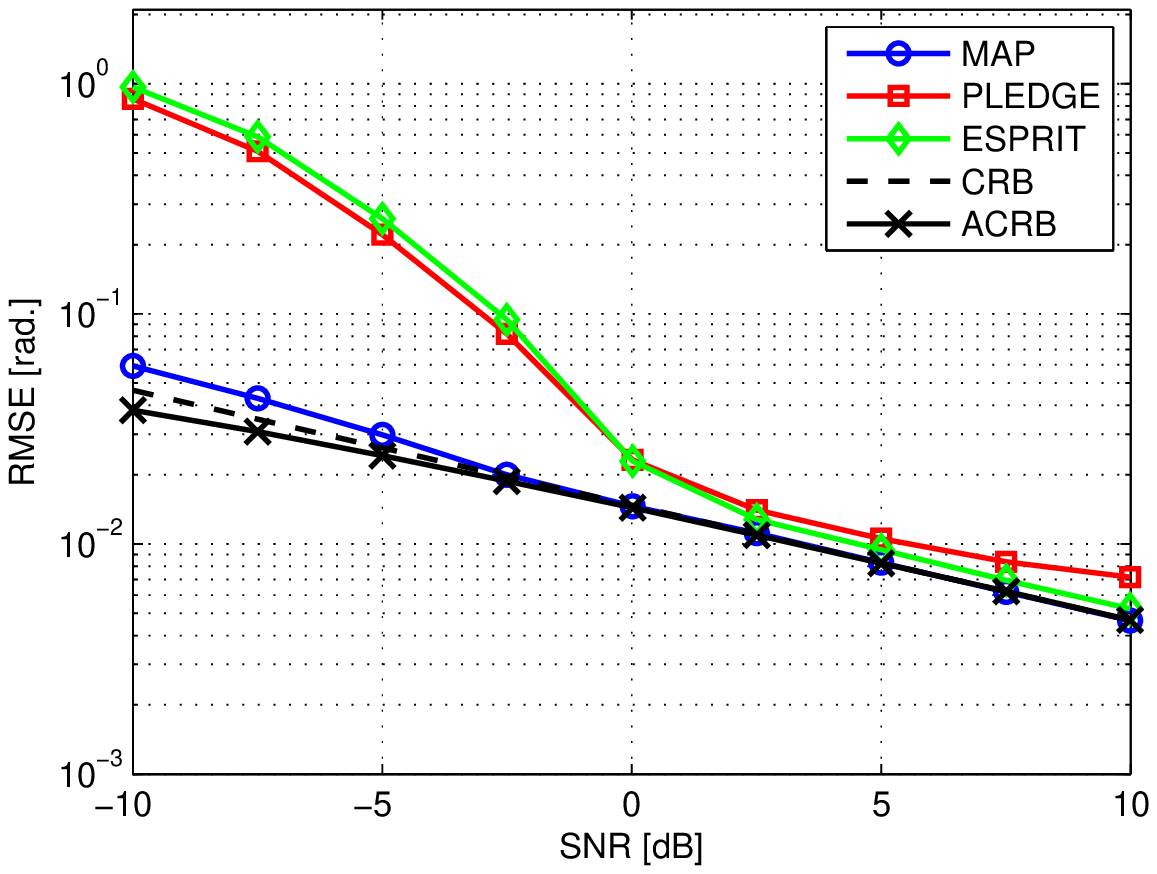}
  \end{center}
  \vspace{-0.5cm}
  \caption{RMSE$(\what{\omega}_2)$ vs SNR for $m=32$.}
  \label{fig:SNR_02_m32}
\end{figure}

\begin{figure}
  \begin{center}
    \includegraphics[width=1.0\columnwidth]{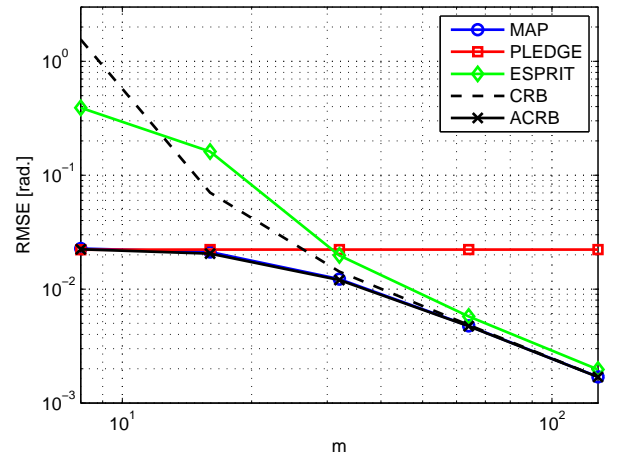}
  \end{center}
  \vspace{-0.5cm}
  \caption{RMSE$(\what{\omega}_1)$ vs $m$ for SNR = 0~dB.}
  \label{fig:M_01_SNR0}
\end{figure}

\begin{figure}
  \begin{center}
    \includegraphics[width=1.0\columnwidth]{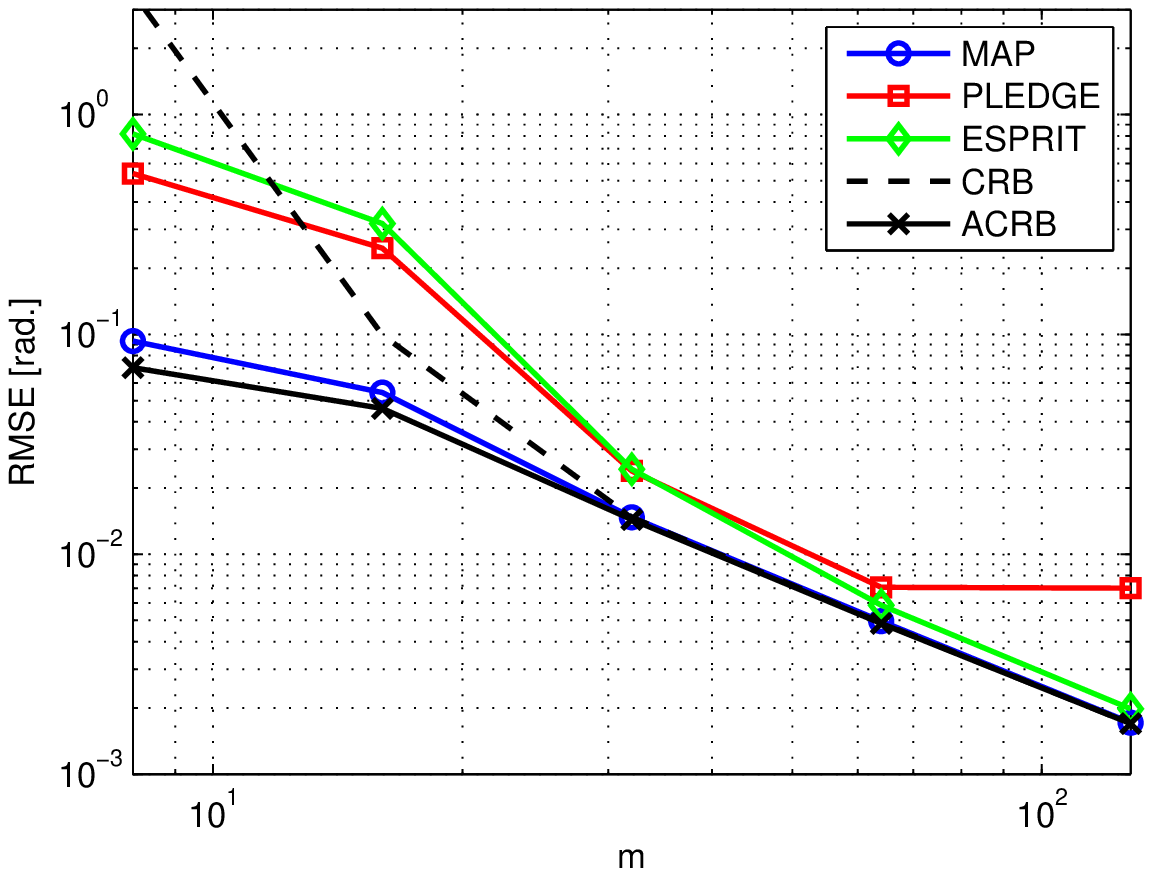}
  \end{center}
  \vspace{-0.5cm}
  \caption{RMSE$(\what{\omega}_2)$ vs $m$ for SNR = 0~dB.}
  \label{fig:M_02_SNR0}
\end{figure}

\begin{figure}
  \begin{center}
    \includegraphics[width=1.0\columnwidth]{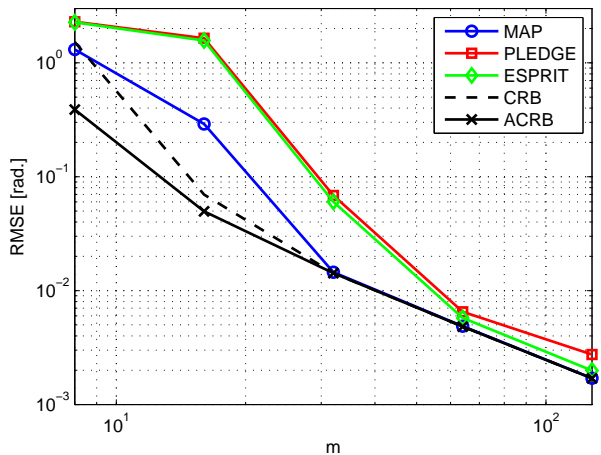}
  \end{center}
  \vspace{-0.5cm}
  \caption{RMSE$(\what{\omega}_3)$ vs $m$ for SNR = 0~dB.}
  \label{fig:M_03_SNR0}
\end{figure}

Next, the number of samples is varied so that $m \in \{ 8, 16, 32,
64, 128 \}$ while fixing SNR=0~dB. The results are displayed in
Figs.~\ref{fig:M_01_SNR0}, \ref{fig:M_02_SNR0} and
\ref{fig:M_03_SNR0}. Again, for $\omega_1$ and $\omega_3$ MAP is able to improve on the prior knowledge; it converges to the CRB from below and follows the ACRB closely. The convergence is also faster for $\omega_3$. The differences between \textsc{Esprit} and \textsc{Pledge} are more visible. The latter
exhibits a smaller gain for $\omega_2$ at low $m$, but is
significantly impaired at high $m$ due to the deterministic modeling of prior knowledge $\omega_1 = \mu_1$.

\emph{Reproducible research:} Code for reproducing empirical results is available at \url{www.ee.kth.se/~davez/rr-line}.

\section{Conclusion}

We have derived the MAP line spectrum estimator in which frequencies are modeled probabilistically. Using the circular von Mises distribution allows for appropriately parameterizing the entire range of uncertainty of the prior knowledge, from complete ignorance to virtual certainty of each frequency. An efficient alternating projections-based solution of the resulting optimization problem was used. The average performance of MAP was then compared to the \textsc{Esprit} and Markov-based \textsc{Pledge} estimators and the Cramér-Rao bound, where its ability to improve on the prior knowledge was demonstrated. MAP would be particularly useful in scenarios with low SNR and/or when few samples are available.

\bibliographystyle{ieeetr}
\bibliography{refs_line}

\end{document}